\documentclass{amsart}
\usepackage{amssymb}
\usepackage{hyperref}
\usepackage{longtable}
\usepackage{orcidlink}
\newcommand{\la}{\lambda}

\newcommand{\Om}{\Omega}
\newcommand{\om}{\omega}

\newcommand{\be}{\begin{equation}}
	\newcommand{\ee}{\end{equation}}
\newcommand{\bea}{\begin{eqnarray}}
	\newcommand{\eea}{\end{eqnarray}}

\newcommand{\bee}{\begin{eqnarray*}}
	\newcommand{\eee}{\end{eqnarray*}}
\newcommand{\ba}{\begin{aligned}}
	\newcommand{\ea}{\end{aligned}}
\newcommand{\bp}{\begin{proof}}
	\newcommand{\ep}{\end{proof}}  
\newcommand{\br}{\begin{remark}}
	\newcommand{\er}{\end{remark}}  
\newcommand{\lb}{\label}

\newtheorem{thm}{Theorem}[section]

\newtheorem{lem}[thm]{Lemma}
\newtheorem{prop}[thm]{Proposition}
\theoremstyle{definition}
\newtheorem{defn}[thm]{Definition}

\theoremstyle{remark}

\usepackage{enumitem}

\newcommand{\Field}{\mathbb{F}}

\begin{document}
	\title[]{ $r$-primitive $k$-normal polynomials over finite fields with last two coefficients prescribed}
	\author[ K. Chatterjee, R. K. Sharma, S. K. Tiwari]{ Kaustav Chatterjee $^*$ \orcidlink{0000-0002-5504-310X}, Rajendra. K. Sharma \orcidlink{0000-0001-5666-4103}, Shailesh. K. Tiwari\orcidlink{0000-0002-7686-1435}}
	\address{Kaustav Chatterjee, Department of Mathematics,  Indian Institute of Technology Patna, 801106, BIHAR, INDIA.}
	\email{kaustav0004@gmail.com}
	\address{R. K. Sharma, Department of Mathematics, South Asian University, New Delhi, Delhi 110068, INDIA.}
	\email{rksharmaiitd@gmail.com}
	\address{Shailesh Kumar Tiwari, Department of Mathematics,  Indian Institute of Technology Patna, 801106, BIHAR, INDIA.}
	\email{shaileshiitd84@gmail.com \& sktiwari@iitp.ac.in}
	\thanks{*email:kaustav0004@gmail.com}
	\maketitle
	\begin{abstract}
		Let $\xi\in\Field_{q^m}$ be an $r$-primitive $k$-normal element over $\Field_q$, where $q$ is a prime power and $m$ is a positive integer. The minimal polynomial of $\xi$ is referred to be the $r$-primitive $k$-normal polynomial of $\xi$ over $\Field_q$. In this article, we study the existence of an $r$-primitive $k$-normal polynomial over $\Field_q$ such that the last two coefficients are prescribed. In this context, first, we prove a sufficient condition which guarantees the existence of such a polynomial. Further, we compute all possible exceptional pairs $(q,m)$ in case of $3$-primitive $1$-normal polynomials for $m\geq 7$.
	\end{abstract}
	\section{Introduction}\lb{S1}
	Let $\Field_{q^m}$ be the $m$-th degree extension of the finite field $\Field_{q}$, where $q$ is a prime power and $m$ is a positive integer. For any positive integer $r$ such that $r|q^m-1$, an element $\xi\in\Field_{q^m}^*$ is said to be an $r$-primitive element if its multiplicative order is $\frac{q^m-1}{r}$. To be specific, $1$-primitive elements are exactly the same as the primitive elements in $\Field_{q^m}$. Also, for $0\leq k\leq m-1$ , any element $\xi\in\Field_{q^m}$ is said to be $k$-normal element over $\Field_q$, if the linear span of the conjugates $\xi$, $\xi^q$, $\xi^{q^2}$,\ldots, $\xi^{q^{m-1}}$ forms a subspace of dimension $m-k$ over $\Field_{q}$. In particular, $0$-normal elements are exactly the same as the normal elements over $\Field_{q}$. An element $\xi$ is said to be $r$-primitive $k$-normal if it is both $r$-primitive and $k$-normal, whereas the corresponding minimal polynomials of an $r$-primitive and $k$-normal element is defined as $r$-primitive and $k$-normal polynomial.
	
	The existence of primitive, normal, or primitive normal elements with specific properties such as predetermined norm and trace over finite fields is a profound research area and a large scale of investigation has been carried out throughout recent years. These studies have been motivated by their applications in fields like cryptography, signal processing, and coding theory. For more information, interested readers can refer to \cite{RH}. The norms and traces of these elements are linked to specific coefficients in their minimal polynomials, which has inspired researchers to investigate special polynomials with predetermined coefficients \cite{HWR, CW, SS, SDC2000, GL2019, AMS2022, MASI22, KHAS}.
	
	In $1992$, Hansen and Mullen conjectured that for any $m\geq 2$, there exists a primitive polynomial of degree $m$, denoted
	by $\mathfrak{P}(x)=x^m-\sigma_{1}x^{m-1}+\ldots+(-1)^m\sigma_{m}$, where the coefficients $\sigma_{n}$ belong to the finite field $\Field_{q}$, except for the cases: $(q,m,n,a) =(q,2,1,0), (4,3,2,0), (4,3,1,0),$ and $(2,4,2,1)$. In $2004$, Fan and Han \cite{SW2004} demonstrated that the Hansen-Mullen conjecture holds approximately by introducing the $p$-adic method, where the only possible exceptions occur when $n=\frac{m+1}{2}$ if $m$ is odd, and when $n=\frac{m}{2}$ or $\frac{m}{2}+1$ if $m$ is even.  In \cite{FHAAECC}, the authors proved that for $m\geq 7$, the conjecture holds true in $\Field_{q^m}$, for even prime power $q$ and  odd positive integer $m$. Later, Cohen \cite{SDC2006} demonstrated that Hansen-Mullen conjecture holds true for all $m\geq 9$, followed by his joint work with Presern \cite{SM2008}, where the authors provided a complete result for all $m\geq 2$. Moreover, in \cite{FHglasgow}, Fan and Han proved that for any degree $m$, there exists a constant $C(m)$ such that we shall always get a primitive polynomial with the first $n=\lfloor\frac{m-1}{2}\rfloor$ coefficients are being prescribed, whenever $q>C(m)$. Numerous researchers have also extended these results to primitive normal polynomials in several studies, such as \cite{FHF2007, FW2009, SF2009, FWDM2009, LG1994}.
	
	 In this article, we extend the result of Fan and Wang \cite{FWDM2009}, where the authors established the existence of a primitive normal polynomial over $\Field_q$ with last two coefficients are prescribed. Here we try to explore all $r$-primitive $k$-normal polynomials with last two coefficients are predetermined.

	This article is structured in the following manner. Definitions, notations, and characteristic functions which are crucial and to be used throughout the paper are mentioned in Section $\ref{S2}$. Further, a sufficient condition followed by a prime sieve is provided to ensure the existence of such $r$-primitive $k$-normal polynomials in Section $\ref{S3}$. As an illustration, we aim to identify finite fields $\Field_{q^m}$ that contain a $3$-primitive $1$-normal polynomial over $\Field_q$ with the last two coefficients prescribed. In Section \ref{S4}, we try to find out the exceptional pairs $(q,m)$ for which $\Field_{q^{m}}$ may not have a  $3$-primitive $1$-normal polynomial over $\Field_{q}$, for $m \geq 7$. When $q<8$ and $m\geq 7$, we show that such a pair will definitely exist apart from $18$ possible choices, whereas for $q\geq 8$ and $ m\geq 8$, the number of such pairs is $20$. Finally, for $m=7$ and $q\geq 8$, an exhaustive computation yields $40$ prime powers as possible exceptions. Moreover, for $m=1,2,3$, we demonstrate that such pairs cannot exist, and for $m=4$, we provide a necessary condition for their possible existence. 
	\setcounter{section}{1}
	\section{Prerequisite}\lb{S2}
	In this section, we recall some basic definitions, lemmas, estimates, and characteristic functions which will be utilized throughout the article. Here $q$ is a prime power and $m $ is a positive integer.
	
	Let $l|q^m-1$ and $\xi\in\Field_{q^m}^*$. Then, we call $\xi$ is $l$-free if and only if gcd$(l,\frac{q^m-1}{ord(\xi)})=1$, where ord$(\xi)$ stands for the order of $\xi$ in the multiplicative group $\Field_{q^m}^*$. Clearly, we have $\xi$ is primitive in $\Field_{q^{m}}$ if and only if it is $(q^m-1)$-free.
	
	Let $f(x)=\sum_{i=0}^{n}{a_{i}x^{i}}\in \Field_{q}[x]$ be such that $a_n\neq 0$ and $\xi\in\Field_{q^m}$. Then, the additive group $\Field_{q^m}$ becomes an $\Field_{q}[x]$-module by the rule 
	$f\circ\xi=\sum_{i=0}^{n}{a_{i}{\xi}^{q^i}}$. Note that, $(x^m-1)\circ\xi=0$ for all $\xi\in\Field_{q^m}$, leading us to define the following. The $\Field_{q}$-order of $\alpha\in\Field_{q^m}$, denoted by Ord$_q(\xi)$, is the monic $\Field_{q}$-divisor $f$ of $x^m-1$ of minimal degree such that $f\circ\xi=0$.
	
	Like as $l$-free elements, any $\xi\in\Field_{q^m}$ is said to be $f$-free if gcd$(f,\frac{x^m-1}{Ord_q(\xi)})=1$. Clearly, $\xi\in\Field_{q^m}$ is normal if and only of it is $(x^m-1)$-free. Following [\cite{SGDD},Theorem 3.2], an element $\xi\in\Field_{q^m}$ is $k$-normal over $\Field_{q}$ if and only if it has $\Field_{q}$-order having degree $m-k$. Using the following lemma, we can easily construct a $k$-normal element with the help of a normal element.
	\begin{lem}\textbf{(\cite{RL}, Lemma 3.1)}\lb{L2.1}
		Let $g\in\Field_{q}[x]$ be a polynomial of degree $k$ such that $g|x^m-1$. Then, for any normal element $\beta\in\Field_{q^m}$ over $\Field_{q}$, we have $\xi=g\circ \beta$ is $k$-normal.
	\end{lem}
	
	Now, let us recall the characters of an abelian group. Let $\mathfrak{A}$ be a finite abelian group and a character $\la$ is a homomorphism from $\mathfrak{A}$ into the set of complex numbers on the unit circle, i.e., $\psi(a_{1}a_{2})=\psi(a_{1})\psi(a_{2})$ for all $a_{1},a_{2}\in \mathfrak{A}$. The character $\psi_{1}$ defined by $\psi_{1}(a)=1$ for all $a\in \mathfrak{A}$, is said to be the trivial character of $A$. Moreover, the collection of all characters of $\mathfrak{A}$, denoted as $\widehat{\mathfrak{A}}$, forms a group under multiplication and $\mathfrak{A}\cong \widehat{\mathfrak{A}}$. 
	
	\begin{lem}\textbf{(\cite{RH}, Theorem 5.4)}
		Let $\psi$ be a nontrivial character of a finite abelain group $\mathfrak{A}$ and $a$ be any nontrivial element in $\mathfrak{A}$. Then $\sum_{a\in\mathfrak{A}}^{}$$\psi(a)=0$ and $\sum_{\psi\in\widehat{\mathfrak{A}}}^{}$$\psi(a)=0$.
	\end{lem}
	We know that the additive character $\la_0$, defined by $\la_0(\xi)=e^{{2\pi i \mathrm{Tr}(\xi)}/p}$, $\mathrm{Tr}$ being the absolute trace function from $\Field_{q^{m}}$ to $\Field_p$, is termed as the \textit{canonical additive character} of $\Field_{q^m}$. Further, for any additive character $\la\in\hat{\Field}_{q^m}$, there exists some $\theta\in\Field_{q^m}$ such that $\la(x)=\la_0(\theta x)$, for all $x\in\Field_{q^m}$. Again, for $\la\in\hat{\Field}_{q^m}$, $g(x)\in\Field_{q}[x]$ and $\xi\in\Field_{q^m}$, $\Field_{q^m}$ forms an $\Field_{q}[x]$-module under the action, $\la\circ g(\xi)=\la(g\circ \xi)$. For $\la\in\hat{\Field}_{q^m}$, the least degree monic divisor $g$ of $x^m-1$ for which $\la\circ g$ is the trivial additive character is said to be the $\Field_{q}$-order of $\la$, denoted by Ord$_q(\la)$. There are exactly $\Phi_q(g)$ characters of $\Field_{q}$-order $g(x)$, where $\Phi_q(g)=|(\Field_{q}[x]/g)^*|$. Further, $\sum_{h|g}^{}\Phi_q(h)=q^{deg(g)}$.

	The following result is given by Weil \cite{AW1948} as described in \cite{TC2006} at $(1.1)$ to $(1.3)$.
	\begin{lem}\lb{L1}
		Assume that $f(x)=\prod_{i=1}^{s}f_{i}(x)^{n_{i}}\in\Field_{q^m}(x)$, where $f_{i}$'s are polynomials over $\Field_{q^m}$ and $n_{i}\in\mathbb{Z}\smallsetminus\{0\}$ such that $f(x)$ that is not of the form of $r(x)^{d}$, where $r(x)\in\Field_{q^m}(x)$ and $d|q^m-1$. If $\psi\in\hat{\Field}_{q^m}^*$ be a multiplicative character of order $d$ and $\mathfrak{D}$ be the sum of the degrees of the $f_i$'s, then $$\Bigg|\sum_{\xi\in\Field_{q^m},f(\xi)\neq \infty }^{}\psi(f(\xi))\Bigg|\leq(\mathfrak{D}-1)q^{m/2}.$$	
	\end{lem}
	\begin{lem}\textbf{{[{\cite{CM2000}}]}}\lb{L2}
		Let $\psi$ be non-trivial multiplicative of order $r$ and $\la$ be a non-trivial additive character on $\Field_{q^n}$. Let $f$ and $g$ be the rational functions in $\Field_{q^m}(x)$ such that  $f\neq yh^r$, for any $y\in\Field_{q^m}$, and $g\neq h^p-h+\beta$, where $h\in\Field_{q^m}(x)$ and $\beta\in\Field_{q^m}$.   Then, we have
		\begin{equation} \nonumber
			\begin{aligned}
				\Bigg|\sum_{\xi\in \Field_{q^{m}}\smallsetminus\mathcal{S}}\psi(f(\xi))\la(g(\xi))\Bigg|\leq(deg(g_\infty)+l+l'-l''-2)q^{m/2},
			\end{aligned}
		\end{equation}
        where $\mathcal{S}$ is set the poles of $f$ and $g$, $(g)_\infty$ is the pole divisor of $g$,  $l$ is the number of distinct zeros and (non-infinite) poles of $f$, $l'$ is the number of distinct poles of $g$ (including $\infty$) and $l''$ is the number of finite poles of $f$ that are poles or zeros of $g$.
	\end{lem}
	
	We now recall the characteristic functions of $e$-free elements and $f$-free elements, where $e|q^m-1$ and $f|x^m-1$. We have the characteristic function determining $e$-free elements of $\Field_{q^m}^*$ is as follows: 
	\begin{equation}\lb{E1}
		\begin{aligned}
			\rho_{e}:\Field_{q^m}^*\mapsto \{0,1\};\xi \mapsto{}\theta(e)\sum_{d|e}\frac{\mu(d)}{\phi(d)}\sum_{\psi_{d}}\psi_{d}(\xi),
		\end{aligned}
	\end{equation}
	where $\theta(e):=\frac{\phi(e)}{e}$, $\psi_{d}$ represents a multiplicative character of order $d$ in $\hat{\Field}_{q^m}^*$ and $\mu$ is the M\"{o}bius function. Further, the characteristic function determining $f$-free elements in $\Field_{q^m}$ is as follows:
	
	\begin{equation}\lb{E2}
		\begin{aligned}
			\kappa_{f}:\Field_{q^m}\mapsto \{0,1\};\xi\mapsto{}\Theta(f)\sum_{h|g}\frac{\mu_{q}(h)}{\Phi_{q}(h)}\sum_{\la_{h}}\la_{h}(\xi),
		\end{aligned}
	\end{equation}
	where $\Theta(g):=\frac{\Phi_{q}(g)}{q^{deg(g)}}$, $\la_{h}$ stands for any additive character of $\Field_{q}$-order $h$ in $\hat{\Field}_{q^m}$ and $\mu_{q}$ is the M\"{o}bius function for the set of polynomials over $\Field_{q}$ is defined as follows: $\mu_{q}(h)=(-1)^r$; if $h$ is the product of $r$ distinct monic irreducible polynomials over $\Field_{q}$, and $0$ otherwise.

	Next to these, there are some characters such as a set of elements with prescribed trace and norm, given by,
	\begin{equation}\lb{E3}
		\tau_{a}:\Field_{q^m}\mapsto \{0,1\};\xi \mapsto \frac{1}{q}\sum_{\la\in{\hat{\Field}_{q}}}^{}\la(\mathrm{Tr}_{\Field_{q^{m}}/\Field_{q}}(\xi)-a),
	\end{equation}
	$$\text{and}$$
	\begin{equation}\lb{E4}
		\begin{aligned}
			\eta_{c}:\Field_{q^m}^*\mapsto \{0,1\};\xi \mapsto \frac{1}{q-1}\sum_{\psi\in{\hat{\Field}_{q}^*}}^{}\psi(\mathrm{N}_{\Field_{q^{m}}/\Field_{q}}(\xi)c^{-1}).
		\end{aligned}
	\end{equation}
	We also need the following definition which can be deduced by the character orthogonality property.
	\begin{defn}
		For any $\xi\in\Field_{q^m}$, we define the following character sum:\\
		\begin{equation}\nonumber
			\begin{aligned}
				I_{0}(\xi)=\frac{1}{q^m}\sum_{\la\in \hat{\Field}_{q^m}}^{}\la(\xi).
			\end{aligned}
		\end{equation}
		
		Observe that $I_{0}(\xi)=1$ if $\xi=0$ and $I_{0}(\xi)=0$ if $\xi\neq0$.
	\end{defn} 
	The following lemma gives us a certain sum, that will be needed in proving our sufficient condition.
	\begin{lem}\textbf{({\cite{JCV2023}}, Lemma 2.5)}\lb{L4}
		Let $g\in\Field_{q}[x]$ be a divisor of $x^m-1$ of degree $k$ and let $\la$ and $\Om$ be additive character. Then\\
		$$
		\frac{1}{q^m}\sum_{\beta\in\Field_{q^m}^{}}\la(\beta)\Om(g \circ \beta)^{-1}=\begin{cases}
			
			1 ~\text{if} ~\la= g\circ \Om,  \\
			0 ~\text{if} ~\la\neq g\circ \Om.
		\end{cases}\\$$
		Here $g\circ\Om(\xi)=\Om(g\circ\xi)$ for all $\xi\in\Field_{q^m}$. Furthermore, for a given additive character $\la$, the set $\hat{g}^{-1}(\la)=\{\la\in\hat{\Field}_{q^m}:\la=g\circ\Om\}$ has $q^{k}$ elements if $Ord(\la)|\frac{x^m-1}{g}$, and it is empty if $Ord(\la)\nmid\frac{x^m-1}{g}$.
	\end{lem}
	\section{General Results}\lb{S3}
	\subsection{Sufficient Condition}
	In this subsection, our goal is to establish a sufficient condition that guarantees the existence of an $r$-primitive $k$-normal polynomial $\mathfrak{P}(x) =x^m-a_{1}x^{m-1}+\ldots+(-1)^{m-1}a_{m-1}x +(-1)^ma_m\in \Field_q[x]$ with the prescribed coefficients $a_{m-1}$ and $a_m$ for $m>1$. Recall that an $r$-primitive $k$-normal element may not have a minimal polynomial of degree $m$ always, first, we will try to establish a sufficient condition that ensures the same.
	The conjugates $\xi$, $\xi^q$, $\xi^{q^2}$,$\ldots$, $\xi^{q^{m-1}}$ of $\xi \in \mathbb{F}_{q^m}$ over $\mathbb{F}_q$ are distinct if and only if the minimal polynomial of $\xi$ has degree $m$. Otherwise, the minimal polynomial will be of degree $d$ a proper divisor of $m$, whereas the conjugates $\xi$, $\xi^q$, $\xi^{q^2}$,$\ldots,\xi^{q^{d-1}}$ repeat $m/d$ times. Now, for any $k$-normal element $\xi\in\Field_{q^m}$ over $\Field_q$, the degree $d$ of the minimal polynomial of $\xi$ can never be less than to $m-k$, and if $m/2<m-k$, then $d$ can never divide $m$. Therefore, if $k<m/2$ gives us $d=m$, that is, $m$ is the degree of  the minimal polynomial of $\xi$.
	Thus, we assume throughout the article that $k<m/2$ ensuring the degree of the minimal polynomial of $\xi$ is exactly $m$. Furthermore, for an irreducible polynomial of degree $m$, we present the following result, which expresses the coefficient $a_{m-1}$ in terms of the trace and norm for $m>1$.  
	\begin{lem}
		Let $m>1$ and $\mathfrak{P}(x)
		=x^m-a_{1}x^{m-1}+\ldots+(-1)^{m-1}a_{m-1}x+(-1)^{m}a_m$ be an irreducible polynomial over $\Field_q$, $\xi$ be a root of $\mathfrak{P}(x)$ in $\Field_{q^m}^*$, where $q$ is prime power. Then $a_{m-1}=\mathrm{N}_{\Field_{q^m}/\Field_q}(\xi)\cdot\mathrm{Tr}_{\Field_{q^m}/\Field_q}(\xi^{-1})$, where $\mathrm{Tr}_{\Field_{q^m}/\Field_q}$ and $\mathrm{N}_{\Field_{q^m}/\Field_q}$ are respectively trace and norm function from $\Field_{q^m}$ to $\Field_q$.
	\end{lem}
	The following lemmas will be useful to prove our main result.
	
	\begin{lem}\textbf{({\cite{AMS2022}, Lemma 3.1})}\lb{L3.2}
		Let $e|q^m-1$, $\delta$=gcd$(e,q-1)$, and $Q_e$ be the largest divisor of $e$ such that gcd$(Q_e,\delta)=1$. Then $\xi\in\Field_{q^m}^*$ is $e$-free if and only if $\xi$ is $Q_e$-free and $\mathrm{N}_{\Field_{q^m}/\Field_{q}}(\xi)$ is $\delta$-free in $\Field_{q}^*$.
	\end{lem} 
	\begin{lem}\textbf{({\cite{MASA 2023}, Lemma 3.1})}\lb{L3.3}
		An element $\xi\in\Field_{q^{m}}$ is $r$-primitive if and only if $\xi=\epsilon^r$ for some primitive $\epsilon\in\Field_{q^{m}}.$ 
	\end{lem}
	\begin{lem}\textbf{({\cite{MASA 2023}, Lemma 3.2})}\lb{L3.4}
		Let $r|q^m-1$ and $w=$gcd$(r,q-1)$. If $\xi$ is an $r$-primitive element, then $\mathrm{N}_{\Field_{q^m}/\Field_{q}}(\xi)$ is $w$-free element $b$ in $\Field_q^*$. Conversely, for any $w$-free element $b\in\Field_{q}^*$, there exists a $r$-primitive element $\xi\in\Field_{q^m}^*$ such that $\mathrm{N}_{\Field_{q^m}/\Field_{q}}(\xi)=b$.
	\end{lem}
	Moreover, we know that $a_m=\mathrm{N}_{\Field_{q^m}/\Field_{q}}(\xi)$. Thus, whenever $k<m/2$, the question of existence of a $r$-primitive $k$-normal polynomial $\mathfrak{P}(x)$ such that $\sigma_{m-1}=a$ and $\sigma_{m}=b$ for any $a,b\in\Field_{q}$ turns into the problem of the existence of a $r$-primitive $k$-normal element $\xi\in\Field_{q^{m}}$ over $\Field_{q}$ such that $\mathrm{Tr}_{\Field_{q^m}/\Field_{q}}(\xi^{-1})=ab^{-1}$ and $\mathrm{N}_{\Field_{q^m}/\Field_{q}}(\xi)=b$, for any $a\in\Field_{q}$, $b\in\Field_{q}^*$. Following lemmas play important role in completion of our main theorem.
    \begin{lem}\lb{L305}
        Let $p,s,q,r,n\in\mathbb{N}$ be such that $q=p^s$ and $r|q^n-1$. Then, for $(\theta,t)\in\Field_{q^n}^*\times\Field_q$, we have $\theta x^r+tx^{-r}\neq h^p-h+c$ for any $h\in\Field_{q^m}(x)$ and $c\in\Field_{q^m}$.   
    \end{lem}
    \bp
If possible, let $\theta x^r+tx^{-r}=h^p-h+c$ for some $h\in\Field_{q^m}(x)$ and $c\in\Field_{q^m}$. Let us write $h=\frac{h_1}{h_2}$, where $h_1$ and $h_2$ are co-prime polynomials over $\Field_{q^m}$. Without loss of generality, suppose that $h_2$ is monic. Here the following possible cases occur: $t=0$ and $t\neq 0$.\\
\textbf{Case 1:} If $t\neq 0$, then we have
\begin{equation}\lb{Eq305}
    (\theta x^{2r}+t)h_{2}^p=x^r(h_{1}^p-h_1h_{2}^{p-1}+ch_{2}^p).
\end{equation}
Since gcd$(x^r,x^{2r}+t)=1$ and gcd$(h_{1}^p-h_1h_{2}^{p-1}+ch_{2}^p,h_{2}^p)=1$, from Equation $(\ref{Eq305})$, we get that $h_2^p=x^r$. This implies that $h_2=x^d$, for some $d\geq 2$, which further implies $r=pd$, a contradiction.  
 \\
\textbf{Case 2:}
If $t=0$, then $h_2$ is constant and thus we get $\theta x^{r}=(h_{1}^p-h_1+c)$, $h_1\in\Field_{q^m}[x]$ and $c\in\Field_{q^m}$. By equating degrees in both sides we get $p|r$, a contradiction. \\

Hence the proof is complete.
\ep
\begin{lem}\lb{L306}
    Let $p,s,q,r,n\in\mathbb{N}$ be such that $q=p^s$ and $r|q^n-1$. Then, for $t\in\Field_{q}^*$, we have $tx^{-r}\neq h^p-h+c$ for any $h\in\Field_{q^m}(x)$ and $c\in\Field_{q^m}$.  
\end{lem}
\bp
The idea of the proof is followed from Lemma \ref{L305} and thus we omit it.
\ep

Assume that $r|q^m-1$ and set $w=$gcd$(r,q-1)$. Let $g$ be a polynomial over $\Field_q$ of degree $k$ such that $g|x^m-1$. Moreover, choose a $w$‑free element $b\in\Field_q^*$ and an arbitrary $a\in\Field_q$. Our objective is to establish a sufficient condition for the existence of $r$‑primitive $k$‑normal elements $\xi$ in $\Field_{q^m}$ that satisfy $\mathrm{Tr}_{\Field_{q^m}/\Field_q}(\xi^{-1})=ab^{-1}$ and $\mathrm{N}_{\Field_{q^m}/\Field_q}(\xi)=b$. To achieve this, it suffices to find a primitive element $\epsilon$ such that $\epsilon^r$ is $k$‑normal, $\mathrm{N}_{\Field_{q^m}/\Field_q}(\xi)=c$ where $c\in\Field_q$ is primitive with $c^r=b$ (see Lemmas \ref{L3.3}, \ref{L3.4}), and $\mathrm{Tr}_{\Field_{q^m}/\Field_q}(\epsilon^{-r})=ab^{-1}$. In fact, we establish a more general statement, from which our desired result follows immediately.

Now, suppose that $e|q^m-1$, $f|x^m-1$, and set $\delta=$gcd$(e,q-1)$. Let $Q_e$ be the greatest divisor of $e$ for which gcd$(Q_e,\delta)=1$, and fix $c\in\Field_q^*$ a $\delta$-free element. The following theorem gives a sufficient criterion that guarantees the existence of a $Q_{e}$ free element $\epsilon$ such that $\epsilon^r=g\circ \beta$, for some $f$ free element $\beta\in\Field_{q^m}$, together with $\mathrm{Tr}_{\Field_{q^m}/\Field_q}(\epsilon^{-r})=ab^{-1}$ and $\mathrm{N}_{\Field_{q^m}/\Field_q}(\epsilon)=c$.
Denote by $\mathfrak{N}_{r,k,a,b,c}(Q_{e},f)$ the number of elements in $\Field_{q^m}$ that satisfy these conditions. Abbreviating $Q:=Q_{q^m-1}$, we provide a criterion on $(q,m)$ to ensure $\mathfrak{N}_{r,k,a,b,c}(Q,x^m-1)>0$.
	\begin{thm}\lb{T3.4}Let $q,m$ be positive integers such that $q$ be a prime power. Also, let $r$ and $k$ be integers such that $r|q^m-1$, and $k<m/2$ such that there exists polynomial $g(x)\in\Field_{q}[x]$ of degree $k$. Then there exists an $r$-primitive $k$-normal element $\xi\in\Field_{q^m}^*$ satisfying $\mathrm{Tr}_{\Field_{q^m}/\Field_q}(\xi^{-1})=ab^{-1}$ and $\mathrm{N}_{\Field_{q^m}/\Field_q}(\xi)=b$, for any prescribed $a\in\Field_q$ and $b\in\Field_{q}^*$ if
		\begin{equation}\nonumber
			\begin{aligned}
				q^{\frac{m}{2}-k-2}>rW(Q)W\Bigg(\frac{x^m-1}{g}\Bigg).
			\end{aligned}
		\end{equation}
	\end{thm}
	\bp  Using the characteristic functions $\rho_{e}$, $\kappa_f$, $I_0$, $\tau_a$ and $\eta_c$, we get that 
	\begin{equation}\nonumber
		\begin{aligned}
			\mathfrak{N}_{r,k,a,b,c}(Q_{e},f)&=\sum_{\epsilon\in\Field_{q^m}^*}^{}\sum_{\beta\in\Field_{q^m}}^{}\rho_{Q_e}(\epsilon)\kappa_f(\beta)I_0(\epsilon^r-g\circ\beta)\tau_{ab^{-1}}(\epsilon^{-r})\eta_{c}(\epsilon)\\
			&=\mathfrak{C}\sum_{d|Q_e,h|f}^{}\frac{\mu(d)\mu_q(h)}{\phi(d)\Phi_q(h)}\sum_{\psi_{d},\la_h,\Om}^{}\mathfrak{T}(\psi_{d},\la_h,\Om),
		\end{aligned}
	\end{equation}
	where $\mathfrak{C}=\frac{\phi(Q_e)\Phi_q(f)}{Q_e\cdot q^{deg(f)+m+1}\cdot(q-1)}$ and 
	\begin{equation}\nonumber
		\begin{aligned}
			\mathfrak{T}(\psi_{d},\la_h,\Om)&=\sum_{\psi\in{\hat{\Field}_q}^*,\la\in{\hat{\Field}_q}}^{}\sum_{\epsilon\in\Field_{q^m}^*}^{}\sum_{\beta\in\Field_{q^m}}^{}\psi_{d}(\epsilon) \la_h(\beta)\Om(\epsilon^r-g\circ\beta)\\&\times\la(\mathrm{Tr}_{\Field_{q^m}/\Field_q}(\epsilon^{-r})-ab^{-1})\psi(\mathrm{N}_{\Field_{q^m}/\Field_q}(\epsilon)c^{-1}).
		\end{aligned}	
	\end{equation}
	Further, let $\la_0$ be the canonical additive character in $\hat{\Field}_q$ and $\psi_{q-1}$ be a multiplicative character of order $q-1$ in ${\hat{\Field}_q}^*$. Thus,
	\begin{equation}\nonumber
		\begin{aligned}
			\mathfrak{T}(\psi_{d},\la_h,\Om)&=\sum_{i=1}^{q-1}\psi_{q-1}(c^{-i})\sum_{t\in\Field_{q}}^{}\la_0(-tab^{-1})\sum_{\epsilon\in\Field_{q^m}^*}^{}\sum_{\beta\in\Field_{q^m}}^{}\psi_{d}(\epsilon) \la_h(\beta)\\&\times\Om(\epsilon^r-g\circ\beta){\la}_0(t\mathrm{Tr}_{\Field_{q^m}/\Field_q}(\epsilon^{-r}))\psi_{q-1}^i(\mathrm{N}_{\Field_{q^m}/\Field_q}(\epsilon))\\&=\sum_{i=1}^{q-1}\psi_{q-1}(c^{-i})\sum_{t\in\Field_{q}}^{}\la_0(-tab^{-1})\sum_{\epsilon\in\Field_{q^m}^*}^{}\sum_{\beta\in\Field_{q^m}}^{}\psi_{d}(\epsilon) \la_h(\beta)\\&\times\Om(\epsilon^r-g\circ\beta)\tilde{\la}_0(t\epsilon^{-r})\tilde{\psi}(\epsilon^i),
		\end{aligned}	
	\end{equation}
	where $\tilde{\psi}=\psi_{q-1}\circ\mathrm{N}_{\Field_{q^m}/\Field_q} $ is a multiplicative character of order $q-1$ in $\Field_{q^m}^*$ and $\tilde{\la}_0=\la_0\circ\mathrm{Tr}_{\Field_{q^m}/\Field_q}$ represents the canonical additive character of $\Field_{q^m}$. Moreover, we have a multiplicative character $\psi_{q^m-1}$ of order $q^m-1$ such that $\tilde{\psi}=\psi_{q^m-1}^{q^m-1/q-1}$, and $\psi_{d}=\psi_{q^m-1}^l$ for some $0\leq l\leq q^m-2$. Also, there exist $\theta,\theta'\in\Field_{q^m}$ such that $\la_h(\beta)=\tilde{\la}_0(\theta'\beta)$ and $\Om(\epsilon)=\tilde{\la}_0(\theta\epsilon)$. Thus,
	\begin{equation}\nonumber
		\begin{aligned}
			\mathfrak{T}(\psi_{d},\la_h,\Om)&=\sum_{i=1}^{q-1}\psi_{q-1}(c^{-i})\sum_{t\in\Field_{q}}^{}\la_0(-tab^{-1})\sum_{\epsilon\in\Field_{q^m}^*}^{}\psi_{q^m-1}(\epsilon^{l+\frac{q^m-1}{q-1}i})\\ &\times\tilde{\la}_0(\theta\epsilon^r+ t\epsilon^{-r})\sum_{\beta\in\Field_{q^m}}^{}\tilde{\la}_0(\theta'\beta-\theta g\circ\beta).
		\end{aligned}
	\end{equation}
	We now write $\mathfrak{R}_{1}(x)=x^{l+\frac{q^m-1}{q-1}i}$, $\mathfrak{R}_{2}(x)=\theta x^r+t x^{-r}$, $\mathfrak{R}_{3}(x)=\theta' x- \theta g\circ x$ and substituting it in the above, yields
	\begin{equation}\nonumber
		\begin{aligned}
			\mathfrak{T}(\psi_{d},\la_h,\Om)&=\sum_{i=1}^{q-1}\psi_{q-1}(c^{-i})\sum_{t\in\Field_{q}}^{}\la_0(-tab^{-1})\sum_{\epsilon\in\Field_{q^m}^*}^{}\psi_{q^m-1}(\mathfrak{R}_{1}(\epsilon))\tilde{\la}_0(\mathfrak{R}_{2}(\epsilon))\\ &\times\sum_{\beta\in\Field_{q^m}}^{}\tilde{\la}_0(\mathfrak{R}_{3}(\beta)).
		\end{aligned}
	\end{equation}
	Before we move further, first notice that
	\begin{equation}\nonumber
		\begin{aligned}
			\sum_{\beta\in\Field_{q^m}}^{}\tilde{\la}_0(\mathfrak{R}_{3}(\beta))=\sum_{\beta\in\Field_{q^m}}\tilde{\la}_0(\theta'\beta-\theta g\circ\beta)=\sum_{\beta\in\Field_{q^m}}\la_h(\beta)\Om(g\circ\beta)^{-1},
		\end{aligned}
	\end{equation}
	where $\la_h(\beta)=\tilde{\la}_0(\theta'\beta)$, $\Om(\beta)=\tilde{\la}_0(\theta\beta)$, and from Lemma \ref{L4}, we have $\underset{\beta\in\Field_{q^m}}{\sum}\la(\beta)\Om(g \circ \beta)^{-1}=q^m$ for all additive characters $\la_h$ and $\Om$ satisfying $\la_h=g\circ \Om$, and the sum is zero elsewhere. Further, for any fixed $\la_h$, there are $q^k$ such characters $\Om$, for given $h|\frac{x^m-1}{g}$. Thus, we can choose $h$ to be a divisor of $\tilde{f}$ only, where $\tilde{f}$=gcd$\Bigg(f,\frac{x^m-1}{g}\Bigg)$.

	Now, let $\psi_{1}$ and $\la_1$ be the trivial multiplicative  character and trivial additive characters respectively. We now write, 
	\begin{equation}\lb{Eq2}
		\begin{aligned}
			\mathfrak{N}_{r,k,a,b,c}(Q_{l},f)
			&=\mathfrak{C}(T_1+T_2+T_3+T_4),
		\end{aligned}
	\end{equation}
	where $T_{1}=\mathfrak{T}(\psi_{1},\la_1,\la_1)$, $T_2=\underset{d|Q_e}{\sum}\frac{\mu(d)}{\phi(d)}\underset{\Om\neq\la_1}{\underset{\psi_{d}}{\sum}}\mathfrak{T}(\psi_{d},\la_1,\Om)$,\\
\begin{equation}\nonumber
		\begin{aligned}
			T_3&=\underset{(d,h)\neq (1,1)}{\sum_{d|Q_e,h|\tilde{f}}^{}}\frac{\mu(d)\mu_q(h)}{\phi(d)\Phi_q(h)}\sum_{\psi_{d},\la_h}^{}\mathfrak{T}(\psi_{d},\la_h,\la_1),
		\end{aligned}
	\end{equation}
and  
\begin{equation}\nonumber
		\begin{aligned}
			T_4&=\underset{h\neq 1}{\sum_{d|Q_e,h|\tilde{f}}^{}}\frac{\mu(d)\mu_q(h)}{\phi(d)\Phi_q(h)}\underset{\Om\neq \la_1}{\sum_{\psi_{d},\la_h}^{}}\mathfrak{T}(\psi_{d},\la_h,\Om).
		\end{aligned}
	\end{equation}
	For $T_1$, we have $l=0$, $\theta'=0$, $\theta=0$, and thus
	\begin{equation}\nonumber
		\begin{aligned}
			T_{1}=&\mathfrak{T}(\psi_{1},\la_1,\la_1)\\=&q^m\sum_{i=1}^{q-1}\psi_{q-1}(c^{-i})\sum_{t\in\Field_{q}}^{}\la_0(-tab^{-1})\sum_{\epsilon\in\Field_{q^m}^*}^{}\psi_{q^m-1}(\epsilon^{\frac{q^m-1}{q-1}i})\tilde{\la}_0(t\epsilon^{-r})\\&=q^m\{(q^m-1)+T_1'+T_1''\},
		\end{aligned}
	\end{equation}
	where
	\begin{equation}\nonumber
		T_1'=\sum_{i=1}^{q-2}\psi_{q-1}(c^{-i})\sum_{\epsilon\in\Field_{q^m}^*}^{}\psi_{q^m-1}(\epsilon^{\frac{q^m-1}{q-1}i}),
	\end{equation} 
	$$\text{and}$$
	\begin{equation}\nonumber
		T_1''=\sum_{i=1}^{q-1}\psi_{q-1}(c^{-i})\sum_{t\in\Field_{q}^*}^{}\la_0(-tab^{-1})\sum_{\epsilon\in\Field_{q^m}^*}^{}\psi_{q^m-1}(\epsilon^{\frac{q^m-1}{q-1}i})\tilde{\la}_0(t\epsilon^{-r}).
	\end{equation}
	Now, $\psi_{q^m-1}^{\frac{q^m-1}{q-1}i}$ is a nontrivial multiplicative character of $\Field_{q^{m}}^*$ for $1\leq i\leq q-2$ and thus we get $T_1'=0$. Following Lemma \ref{L306}, for $t\in\Field_{q}^*$, $tx^{-r}\neq h(x)^{p}-h(x)+u$ for any $h(x)\in\Field_{q^m}(x)$, $u\in\Field_{q^m}$, and so by Lemma \ref{L2}, we have $|T_1''|\leq r(q-1)^2q^{m/2}$. Combining these two, we get $|T_1-q^m(q^m-1)|\leq r(q-1)^2q^{3m/2}$.
	
	We now try to get an estimate of $T_2$, where, 
	\begin{equation}\nonumber
		\begin{aligned}
			T_2=\sum_{d|Q_e}^{}\frac{\mu(d)}{\phi(d)}\sum_{\psi_{d},\Om\neq\la_1}^{}\mathfrak{T}(\psi_{d},\la_1,\Om).
		\end{aligned}
	\end{equation}
	For this, $\la_h$ is the trivial additive character $\la_1$, that is, $\theta'=0$ and $\Om$ is non-trivial additive character, i.e., $\theta\neq 0$. Therefore, by Lemma \ref{L4}, we have the following
	\begin{equation}\nonumber
		\begin{aligned}
			\sum_{\Om\in\hat{\Field}_{q^m},\Om\neq \la_1}^{}\sum_{\beta\in\Field_{q^m}}\tilde{\la}_0(\theta'\beta-\theta g\circ\beta)=(q^k-1)q^m.
		\end{aligned}
	\end{equation}
	Thus, 
	\begin{equation}\nonumber
		\begin{aligned}
			T_2&=\sum_{d|Q_e}^{}\frac{\mu(d)}{\phi(d)}\sum_{\psi_{d},\Om\neq\la_1}^{}\mathfrak{T}(\psi_{d},\la_1,\Om)\\&=(q^k-1)q^m\sum_{d|Q_e}^{}\frac{\mu(d)}{\phi(d)}\sum_{\psi_{d}}^{}\sum_{i=1}^{q-1}\psi_{q-1}(c^{-i})\sum_{t\in\Field_{q}}^{}\la_0(-tab^{-1})\\&\times\sum_{\epsilon\in\Field_{q^m}^*}^{}\psi_{q^m-1}(\epsilon^{l+\frac{q^m-1}{q-1}i})\tilde{\la}_0(\theta\epsilon^r+ t\epsilon^{-r}).
		\end{aligned}
	\end{equation}
	Since $\theta\neq 0$, by Lemma \ref{L305}, $\theta x^r+t x^{-r}\neq 
	{h}^{p}-h+c$, for $h\in\Field_{q^m}(x)$ and $c\in\Field_{q^m}(x)$. Thus, following Lemma \ref{L2}, we get that $|T_2|\leq r(q^k-1)(q-1)q^{\frac{3m}{2}+1}W(Q_e)$.
	
	We now estimate $T_{3}$. For this, we have 
	\begin{equation}\nonumber
		\begin{aligned}
			T_3&=\underset{(d,h)\neq (1,1)}{\sum_{d|Q_e,h|\tilde{f}}^{}}\frac{\mu(d)\mu_q(h)}{\phi(d)\Phi_q(h)}\sum_{\psi_{d},\la_h}^{}\mathfrak{T}(\psi_{d},\la_h,\la_1)\\&=\underset{(d,h)\neq (1,1)}{\sum_{d|Q_e,h|\tilde{f}}^{}}\frac{\mu(d)\mu_q(h)}{\phi(d)\Phi_q(h)}\sum_{\psi_{d},\la_h}^{}\sum_{i=1}^{q-1}\psi_{q-1}(c^{-i})\sum_{t\in\Field_{q}}^{}\la_0(-tab^{-1})\\&\times\sum_{\epsilon\in\Field_{q^m}^*}^{}\psi_{q^m-1}(\epsilon^{l+\frac{q^m-1}{q-1}i})\tilde{\la}_0(t\epsilon^{-r})\sum_{\beta\in\Field_{q^m}}\tilde{\la}_0(\theta'\beta).
		\end{aligned}
	\end{equation}
	Here comes two possible scenarios, $h\neq 1$, that is, $\theta'\neq 0$ or $h=1$, that is, $\theta'=0$. The first one gives us $T_3=0$. So we consider only $h=1$, that is, $\theta'=0$. Then $d\neq 1$ and thus 
	\begin{equation}\nonumber
		\begin{aligned}
			T_3&=q^m\underset{d\neq 1}{\sum_{d|Q_e}^{}}\frac{\mu(d)}{\phi(d)}\sum_{\psi_{d}}^{}\Bigg\{\sum_{i=1}^{q-1}\psi_{q-1}(c^{-i})\sum_{t\in\Field_{q}}^{}\la_0(-tab^{-1})\\&\times\sum_{\epsilon\in\Field_{q^m}^*}^{}\psi_{q^m-1}(\epsilon^{l+\frac{q^m-1}{q-1}i})\tilde{\la}_0(t\epsilon^{-r})\Bigg\}\\&=q^m(T_3'+T_3''),
		\end{aligned}
	\end{equation}
	where 
	\begin{equation}\nonumber
		\begin{aligned}
			T_3'=\underset{d\neq 1}{\sum_{d|Q_e}^{}}\frac{\mu(d)}{\phi(d)}\sum_{\psi_{d}}^{}\Bigg\{\sum_{i=1}^{q-1}\psi_{q-1}(c^{-i})\sum_{\epsilon\in\Field_{q^m}^*}^{}\psi_{q^m-1}(\epsilon^{l+\frac{q^m-1}{q-1}i})\Bigg\},
		\end{aligned}
	\end{equation}
	$$\text{and}$$
	\begin{equation}\nonumber
		\begin{aligned}
			T_3''&=\underset{d\neq 1}{\sum_{d|Q_e}^{}}\frac{\mu(d)}{\phi(d)}\sum_{\psi_{d}}^{}\Bigg\{\sum_{i=1}^{q-1}\psi_{q-1}(c^{-i})\sum_{t\in\Field_{q}^*}^{}\la_0(-tab^{-1})\\&\times\sum_{\epsilon\in\Field_{q^m}^*}^{}\psi_{q^m-1}(\epsilon^{l+\frac{q^m-1}{q-1}i})\tilde{\la}_0(t\epsilon^{-r})\Bigg\}.
		\end{aligned}
	\end{equation}
	Now, let us first evaluate $T_3'$. For this assume that $x^{l+\frac{q^m-1}{q-1}i}=ch(x)^{q^m-1}$, for $h\in\Field_{q^m}(x)$, $c\in\Field_{q^m}^*$. Then, by comparison of the power's of $x$, we get that $l+\frac{q^m-1}{q-1}i=(q^m-1)t$, for some positive integer $t$. Following \cite{AMS2022}, we get that $l=0$, a contradiction. Therefore, we have $x^{l+\frac{q^m-1}{q-1}i}\neq H(x)^{q^m-1}$, for any $H(x)\in\Field_{q^m}(x)$, and thus $T_3'=0$. Again, for $t\in\Field_{q}^*$, $tx^{-r}\neq {h(x)}^{p}-h(x)+u$, for $h(x)\in\Field_{q^m}(x)$, $u\in\Field_{q^m}$. Thus, by Lemma \ref{L2}, it follows that $|T_3|\leq r(q-1)^2q^\frac{3m}{2}(W(Q_e)-1)$.
	
	Finally, we estimate $T_4$. For this, we use
	\begin{equation}\nonumber
		\begin{aligned}
			T_4&=\underset{h\neq 1}{\sum_{d|Q_e,h|\tilde{f}}^{}}\frac{\mu(d)\mu_q(h)}{\phi(d)\Phi_q(h)}\underset{\Om\neq \la_1}{\sum_{\psi_{d},\la_h}^{}}\mathfrak{T}(\psi_{d},\la_h,\Om)\\&=\underset{h\neq 1}{\sum_{d|Q_e,h|\tilde{f}}^{}}\frac{\mu(d)\mu_q(h)}{\phi(d)\Phi_q(h)}\underset{\Om\neq \la_1}{\sum_{\psi_{d},\la_h}^{}}\sum_{i=1}^{q-1}\psi_{q-1}(c^{-i})\sum_{t\in\Field_{q}}^{}\la_0(-tab^{-1})\\&\times\sum_{\epsilon\in\Field_{q^m}^*}^{}\psi_{q^m-1}(\epsilon^{l+\frac{q^m-1}{q-1}i})\tilde{\la}_0(\theta\epsilon^r+t\epsilon^{-r})\sum_{\beta\in\Field_{q^m}}^{}\tilde{\la}_0(\theta'\beta-\theta g\circ\beta).
		\end{aligned}
	\end{equation}
	By using the same theory as we did before Equation (\ref{Eq2}), we get $\sum_{\beta\in\Field_{q^m}}^{}\tilde{\la}_0(\theta'\beta-\theta g\circ\beta)=q^m$ if $\la_h=g\circ \Om$ otherwise this sum is zero, and this will happen for $q^k$ such $\Om$
	for every $\la_h$. Thus, we get 
	\begin{equation}\nonumber
		\begin{aligned}
			T_4&=q^{m+k}\underset{h\neq1}{\sum_{d|Q_e,h|\tilde{f}}^{}}\frac{\mu(d)\mu_q(h)}{\phi(d)\Phi_q(h)}\underset{\Om\neq\la_1,\la_h=g\circ\Om}{\sum_{\psi_{d},\la_h}^{}}\sum_{i=1}^{q-1}\psi_{q-1}(c^{-i})\\&\times\sum_{t\in\Field_{q}}^{}\la_0(-tab^{-1})\sum_{\epsilon\in\Field_{q^m}^*}^{}\psi_{q^m-1}(\epsilon^{l+\frac{q^m-1}{q-1}i})\tilde{\la}_0(\theta\epsilon^r+t\epsilon^{-r}).
		\end{aligned}
	\end{equation}
	Since $\Om\neq\la_1$, we must have $\theta\neq 0$ and thus it follows that  $\theta x^r+t x^{-r}\neq{h(x)}^{p}-h(x)+\beta$, for $h(x)\in\Field_{q^m}(x)$ and $\beta\in\Field_{q^m}$. By Lemma \ref{L2}, we have
	\begin{equation}\nonumber
		|T_4|\leq r(q-1)q^{\frac{3m}{2}+k+1}W(Q_e)(W(\tilde{f})-1).
	\end{equation} 
	Now, from the definition of $\mathfrak{N}_{r,k,a,b,c}(Q_{e},f)$, we get that
	\begin{equation}\nonumber
		\begin{aligned}
			\mathfrak{N}_{r,k,a,b,c}(Q_{e},&f)\geq \mathfrak{C}\{q^m(q^m-1)- r(q-1)^2q^{\frac{3m}{2}}-r(q^k-1)(q-1)q^{\frac{3m}{2}+1}W(Q_e)\\&-r(q-1)^2q^\frac{3m}{2}(W(Q_e)-1)-r(q-1)q^{\frac{3m}{2}+k+1}W(Q_e)(W(\tilde{f})-1)\}\\&>\mathfrak{C}\{q^{2m}-rq^{\frac{3m}{2}+k+2}W(Q_e)W(\tilde{f})\}.
		\end{aligned}
	\end{equation}
Thus, we get the sufficient condition for $\mathfrak{N}_{r,k,a,b,c}(Q_{e},f)$ to be positive is given by 
	\begin{equation}\lb{E6}
		q^{\frac{m}{2}-2-k}>rW(Q_e)W(\tilde{f}).
	\end{equation}Now, in particular, if we choose $e=q^m-1$ and $f=x^m-1$, then $\tilde{f}=\frac{x^m-1}{g}$ and the sufficient condition for the existence of an $r$-primitive $k$-normal element $\xi=\epsilon^r$ such that $\mathrm{Tr}_{\Field_{q^m}/\Field_{q}}(\xi^{-1})=ab^{-1}$ and $\mathrm{N}_{\Field_{q^m}/\Field_{q}}(\xi)=b$, for predetermined $a,b\in\Field_{q}$ with $b\neq 0$, is as follows,
	\begin{equation}\lb{E7}
		\begin{aligned}
			q^{\frac{m}{2}-k-2}>rW(Q)W\Bigg(\frac{x^m-1}{g}\Bigg).
		\end{aligned}
	\end{equation}
	In other words, Inequality (\ref{E7}) serves as a sufficient condition which ensures the existence of an $r$-primitive $k$-normal polynomial $\mathfrak{P}(x)$ of degree $m$ over $\Field_{q}$ such that $\sigma_{m-1}$ and $\sigma_{m}$ are both prescribed.\ep
	\subsection{The prime sieve.} 
	This aim of this section is to further improve the condition provided in the Theorem \ref{T3.4}. In the subsequent two results, the sieving technique that is given, is similar as the previous works of primitive and normal elements. The proofs of these two are same as given in \cite{MAS2022,JCV2023}, and thus we omit the proof here.
	\begin{lem}
		Let $r|q^m-1$ be a positive integer and $k$ be a non-negative integer. Let $e|Q$ and $\{p_1,p_2,\ldots,p_s\}$ be the largest set of prime divisors of $Q$  such that $p_i\nmid e$ for all $1\leq i\leq s_1$. Furthermore, let $f$ be a divisor of $x^m-1$ and $\{f_1\,f_2,\ldots,f_{s_2}\}$ be the largest set of irreducible factors of $x^m-1$ such that $f_i\nmid f$ for all $1\leq i\leq s_2$. Then, we have
		\begin{equation}\nonumber
			\begin{aligned}
				\mathfrak{N}_{r,k,a,b,c}(Q,x^m-1)&\geq \sum_{i=1}^{s_1}\mathfrak{N}_{r,k,a,b,c}(ep_i,f)+\sum_{j=1}^{s_2}\mathfrak{N}_{r,k,a,b,c}(e,ff_j)\\&-(s_1+s_2-1)\mathfrak{N}_{r,k,a,b,c}(e,f).
			\end{aligned}
		\end{equation}
	\end{lem}
	\begin{prop}\lb{P3.6}
		We adopt the notations from the above lemma. Further, we define
		\begin{equation}\nonumber
			\delta=1-\sum_{i=1}^{s_1}\frac{1}{p_i}-\sum_{j=1}^{s_2}\frac{1}{q^{deg(f_j)}}
		\end{equation}
		$$\text{and}$$	
		\begin{equation}\nonumber
			\Delta=\frac{s_1+s_2-1}{\delta}+2.
		\end{equation}
		Then, provided $\delta>0$, $\mathfrak{N}_{r,k,a,b,c}(Q,x^m-1)>0$ whenever
		\begin{equation}\lb{E8}
			\begin{aligned}
				q^{m/2-k-2}>rW(e)W(\tilde{f})\Delta.
			\end{aligned}
		\end{equation} 
	\end{prop}
    
	\section{Existence of $3$-primitive $1$-normal elements with \\prescribed norm and trace}\lb{S4}	
	In this section, we aim to identify those finite fields where we can find $3$-primitive $1$-normal polynomials $\mathfrak{P}(x) =x^m-a_{1}x^{m-1}+\ldots+(-1)^{m-1}a_{m-1}x +(-1)^ma_m\in \Field_q[x]$ such that the coefficients $a_{m-1},a_m$ are prescribed. In other words, we shall find a $3$-primitive $1$-normal element $\xi\in\Field_{q^{n}}$ such that $\mathrm{Tr}_{\Field_{q^{m}}/\Field_{q}}(\xi^{-1})=ab^{-1}$ and $\mathrm{N}_{\Field_{q^{m}}/\Field_{q}}(\xi)=b$, for any $a\in\Field_q$ and $b\in\Field_q^*$. To simplify the matter, we shall denote the set $\Gamma_{a,b}(r,k)$ as the number of pairs $(q,m)$ such that the field $\Field_{q^m}$ contains at least one $r$-primitive $k$-normal element $\xi$ such that $\mathrm{Tr}_{\Field_{q^{m}}/\Field_{q}}(\xi^{-1})=ab^{-1}$ and $\mathrm{N}_{\Field_{q^{m}}/\Field_{q}}(\xi)=b$.

	We now recall that for $x$, a positive integer (or a polynomial over $\Field_{q}$), $W(x)$ stands for the number of square free divisors (or square free factors) of $x$. Before we proceed further, we consider the following two lemmas, where bounds of $W(n)$ and $W(x^m-1)$ are provided, for some positive integer $n$. 
	
	\begin{lem}\textbf{[{\cite{SS}}, Lemma 3.7]}\lb{L4.1}
		Let $r>0$ be a real number and $n$ be a positive integer. Then $W(n)<\mathcal{C}_r\cdot n^{\frac{1}{r}}$, where $\mathcal{C}_r=\frac{2^w}{{({p_{1}p_{2}\ldots p_{w})}}^{\frac{1}{r}}}$ and $p_{1},p_{2},\ldots,p_{w}$ are primes $\leq 2^{r}$ that divide $n$.
	\end{lem}
	\begin{lem}\textbf{[{\cite{HWR}}, Lemma 2.9]}\lb{L4.2}
		Suppose $q$ is a prime power, $m$ is a natural number and $m'$=gcd$(m,q-1)$. Then $W(x^m-1)\leq 2^{\frac{1}{2}\{m+m'\}}$, which gives $W(x^m-1)\leq 2^m$. Further, $W(x^m-1)=2^m$ if and only if $m|q-1$. In addition, if  $m\nmid q-1$, then $W(x^m-1)\leq 2^{\frac{3m}{4}}$.
	\end{lem}
	We now assume that $r=3$, $k=1$ and thus the inequalities (\ref{E7}) and (\ref{E8}) turns into the followings, 
	\begin{equation}\lb{E9}
		\begin{aligned}
			q^{\frac{m}{2}-3}>3W(Q)W\Bigg(\frac{x^m-1}{g}\Bigg),
		\end{aligned}
	\end{equation}
	$$\text{and}$$
	\begin{equation}\lb{E10}
		\begin{aligned}
			q^{\frac{m}{2}-3}>3W(e)W(\tilde{f})\Delta.
		\end{aligned}
	\end{equation} 
	Clearly, $3\nmid q^m-1$ implies that $(q,m)\notin\Gamma_{a,b}(3,1)$ and thus we deal with the pairs $(q,m)$ such that $3|q^m-1$. Although, the above inequalities are not sufficient to deal with the scenarios $m<7$. Henceforth, we consider the cases $m\geq 7$ only, and our goal is to explicitly find the pairs $(q,m)\in\Gamma_{a,b}(3,1)$. Also, note that $W(Q)\leq W(q^m-1)$.
	\subsection{\textbf{Cases} $q<8$.} In the following lemma, we address the cases where $q=4,5,7$ and $m\geq 7$. 
	\begin{lem}
		Let $m\geq 7$ be an integer and $q=4,5,7$ be such that $3|q^m-1$. Then $(q,m)\in\Gamma_{a,b}(3,1)$ unless possibly for the following pairs.
		\begin{itemize}
			\item[(i)] $q=4$ and $m=7, 8, 9, 10, 12$.
			\item[(ii)] $q=5$ and $m=8, 10, 12$.
			\item[(iii)] $q=7$ and $m=7, 8, 9, 12$.
		\end{itemize}
	\end{lem}
	\bp
	For $q=4,5,7$ and $m\geq 7$, we have $m\nmid q-1$ and thus $W(x^m-1)\leq 2^{3m/4}$. Therefore, Inequality (\ref{E9}) holds if,
	\begin{equation}
		\begin{aligned}
			m>\frac{\log(3\mathcal{C}_rq^3)}{(1/2-1/r)\log{q}-(3/4)\log{2}}.
		\end{aligned}
	\end{equation}
	Clearly, when $q=4$, for $r=10.3$, the above inequality holds for $m\geq 757$, when $q=5$, for $r=8.4$, then it holds for $m\geq 172$ and when $q=7$, for $r=8.2$, then the inequality holds for $m\geq 75$. For these pairs, we test the inequality (\ref{E9}) and get $(q,m)\in\Gamma_{a,b}(3,1)$ unless possibly $q=4$ and $m=7,8,9,10,11,12,14,15,18$; $q=5$ and $m=8,10,12,16,24$; $q=7$ and $m=7,8,9,10,12$. Among these, the pairs that satisfy $(\ref{E10})$, are as given by $(4,11)$, $(4,14)$, $(4,15)$, $(4,18)$, $(5,16)$, $(5,24)$, and $(7,10)$. Hence the proof is complete.
	\ep
	Next, we discuss about the possible exceptions when $q=2$. Let $q$ be a prime power and $m=q^i\cdot m'$ be a positive integer such that gcd$(m',q)=1$; $i\geq 0$, and $t$ be the order of $q$ modulo $m'$. Then, $x^{m'}-1$ can be factorized in terms of irreducible polynomials of degree $<t$. Let $\rho(q,m')$ be ratio of the number of the irreducible factors of $x^{m'}-1$ to $m'$. Then $m'\rho(q,m')=m\rho(q,m)$. 
	\begin{lem}\textbf{[\cite{SS}, Lemma $7.1$]}\lb{L4.4}
		If $m>4$ and $2\nmid m$, then  $\rho(2,5)=1/5$, $\rho(2,9)=1/9$, $\rho(2,21)=4/21$ and $\rho(2,m)\leq 1/6$ for $m\neq 5,9,21$. 
	\end{lem} 
	\begin{lem}
		Let $q=2$ and $m\geq 7$ be a positive integer such that $3|2^m-1$. Then, $(2,m)\in\Gamma_{a,b}(3,1)$  unless $m=8,10,12,14,16,18$. 
	\end{lem}
	\bp
	Clearly, $m$ must be even positive integer and express $m$ in the form $m'\cdot 2^i$, for some $i\geq 1$ and $2\nmid m'$. Then $x^m-1$ and $x^{m'}-1$ have same factors over $\Field_{q}$. So, there are two possibilities either $m'|q^2-1$ or $m'\nmid q^2-1$.
	
	Firstly, we assume that $m'|q^2-1$. This gives that $m'=1$ or $3$, that is, $x^{m'}-1=x-1$ or $(x-1)(x^2+x+1)$. Choose $e=Q$ and $f=\frac{x^{m'}-1}{x-1}$. Then by Proposition \ref{P3.6}, we have $\delta=\frac{1}{2}$, $\Delta=2$. Further by using Lemma \ref{L4.1}, we get $(q,m)\in\Gamma_{a,b}(3,1)$, if
	\begin{equation}\nonumber
		\begin{aligned}
			q^{\frac{m}{2}-3}> 12\cdot\mathcal{C}_r \cdot q^{\frac{m}{r}}.
		\end{aligned}
	\end{equation}
	Take $r=3$, which gives that $(q,m)\in\Gamma_{a,b}(3,1)$ if $2^m> (96\cdot 2.7)^6$, which holds for $m\geq 49$. Further, when $7\leq  m\leq 48$, the values of $m$ for which we need to verify are given by, $m=8, 12, 16, 24, 32, 48$. Fore these values, after verifying Inequalities $(\ref{E9})$ and $(\ref{E10})$, it follows that $(2,m)\in\Gamma_{a,b}(3,1)$ unless $m=8,12, 16$. 
	
	Secondly, we assume that $m'\nmid q^2-1$.  Then, we have $m'\geq 5$ and $t>2$. Here, we choose $e=Q$ and $f$ as the product of the irreducible factors of degree $<t$. Then, following [\cite{MASI22}, Lemma $10$], it follows that $\delta>1-1/t>0$ and $\Delta\leq m'\leq m$. Further, note that $\tilde{f}$=gcd$(f,\frac{x^m-1}{g})=f$ and thus $W(\tilde{f})=2^{m'\rho(q,m')}=2^{m\rho(q,m)}$. From Lemma \ref{L4.4}, $m\rho(q,m)\leq m/6$ and hence $(q,m)\in\Gamma_{a,b}(3,1)$, if 
	\begin{equation}\nonumber
		\begin{aligned}
			q^{\frac{m}{2}-3}> 3\mathcal{C}_r q^{\frac{m}{r}}2^{\frac{m}{6}}m.
		\end{aligned}
	\end{equation}
	Choosing $r=6$, the above inequality holds for all $m\geq 99$. For $7\leq m\leq 98$, after verifying Inequality $(\ref{E9})$ and the pairs $(2,m)$ not satisfying the same, we verify Inequality $(\ref{E10})$. Hence, we get that $(2,m)\in\Gamma_{a,b}(3,1)$ unless $m=10,14,18$.
	\ep
	
	\subsection{ The Cases $q\geq 8$} In the following lemma, we focus on the cases $m\geq 7$ and $q\geq 8$.
	\begin{lem}
		Let $m\geq 12$ be a positive integer and $q\geq 8$ be a prime power be such that $3|q^m-1$, then $(q,m)\in\Gamma_{a,b}(3,1)$. 
	\end{lem}
	\bp
	Let $q\geq 8$ and $m\geq 12$. By using Lemmas \ref{L4.1} and \ref{L4.2}, Inequality (\ref{E9}) holds if $q^{m/2-3}>3~\mathcal{C}_rq^{m/r}2^{m-1}$, or equivalently
	\begin{equation}\lb{E12}
		\begin{aligned}
			q>(3\cdot 2^{m-1}\mathcal{C}_r)^{\frac{2r}{(m-6)r-2m}},
		\end{aligned}
	\end{equation}
	which holds whenever $r>\frac{2m}{m-6}$. The pairs $(q,m)$ that satisfy Inequality (\ref{E12}) for certain values of $r$, are given in Table \ref{Table1} or in other words, the pairs $(q,m)\in\Gamma_{a,b}(3,1)$.
	\begin{center}
		\begin{table}[h]
			\centering
			\caption{Pairs $(q,m)\in\Gamma_{a,b}(3,1)$.}
			\setlength{\tabcolsep}{40pt} 
			\begin{tabular}{|l|l|}
				\hline $r$ & $(q,m)$  \\
				\hline 8.5 & $q\geq 8$ and $m\geq 175$  \\
				7.7 & $q\geq 13$ and $m\geq 64$ \\
				7.5 & $q\geq  21$ and $m\geq 39$ \\
				7.5 & $q\geq 84$ and $m\geq 22$  \\
				7.5 & $q\geq 354$ and $m\geq 17$  \\
				7.5 & $q\geq 2747$ and $m\geq 14$  \\
				7.2 & $q\geq 8822$ and $m\geq 13$  \\
				7.2 & $q\geq 61429$ and $m\geq 12$  \\
				\hline
			\end{tabular}
			\label{Table1}
		\end{table}
	\end{center}
	The values of $q$ and $m$ that do not satisfy Inequality $(\ref{E12})$ are given in Table \ref{Table2}. For the remaining pairs in the corresponding table, we enumerate all prime powers $q$ within the specified range. Moreover, for each $m$ within the given range, we first verify Inequality $(\ref{E9})$. If it fails to hold, we then test the Inequality $(\ref{E10})$ by taking certain choices of $e$ and $f$. As a result, we find that all $(q,m)\in\Gamma_{a,b}(3,1)$, for all $n\geq 12$ and $q\geq 8$. 
	\ep
	\begin{center}
		\begin{table}[h]
			\centering
			\caption{Pairs $(q,m)$ that do not satisfy Inequality $(\ref{E12})$.}
			\setlength{\tabcolsep}{40pt} 
			\begin{tabular}{|l|l|}
				\hline $q$ & $m$  \\
				\hline $8\leq q < 13$ & $64\leq m< 175$  \\
				$8\leq q <21$& $39\leq m< 64$\\
				$8\leq q< 84$ & $22\leq m< 39$ \\
				$8\leq q< 354$ & $17\leq m< 22$  \\
				$8\leq q< 2747$ & $14\leq m< 17$  \\
				$8\leq q< 8822$ & $13\leq m< 14$  \\
				$8\leq q< 61429$ & $12\leq m< 13$  \\
				\hline
			\end{tabular}
			\label{Table2}
		\end{table}
	\end{center}
	Now, for $7\leq m\leq 11$, we test Inequality $(\ref{E12})$ and get the following lower bounds of $q$ (which are very large) such that $(q,m)\in\Gamma_{a,b}(3,1)$.
	\begin{itemize}
		\item[1.] When $m=11$ and $1.69\times 10^6$ (choosing $r=7.5$).
		\item[2.] When $m=10$ and $q\geq 1.28\times 10^9$ (choosing $r=7.9$). 
		\item[3.] When $m=9$ and $q\geq 7.23\times 10^{16}$ (choosing $r=8.6$).
		\item[4.] When $m=8$ and $q\geq 1.12\times 10^{58}$ (choosing $r=10.1$).
		\item[5.] When $m=7$ and $q\geq 3.59\times 10^{3177}$ (choosing $r=15.7$).
	\end{itemize}
	For the remaining cases, it is computationally challenging to provide a complete list of all the prime powers. Therefore, in the following lemmas, we start by lowering these bounds. Then, within the reduced range, we list all the prime powers to test the validity of pairs $(q, m) \in \Gamma_{a,b}(3,1)$ using Inequality $(\ref{E9})$ and the sieving inequality $(\ref{E10})$.
	
	Following [\cite{ARS2018} Lemma $5.1$], we have deduced the following lemma, which will help us in finding the exceptions in the remaining cases. 
	\begin{lem}\lb{L4.5}
		Let $M\in\mathbb{N}$ such that $\om(M)\geq 9632$. Then $W(M)<M^{1/15}$.
	\end{lem} 
	From now onwards, we assume that $7\leq m\leq 11$ and $q\geq 8$. Further, let $\om(q^m-1)\geq 9632$ and then, by using Lemma \ref{L4.5} and Inequality $(\ref{E9})$, $(q,m)\in\Gamma_{a,b}(3,1)$ if $q^{\frac{m}{2}-3}>3\cdot 2^{10}\cdot q^{\frac{m}{15}}$, that is, if $q^m>(3\cdot 2^{10})^{\frac{30m}{13m-90}}$. But, $m\geq 7$ gives us $\frac{30m}{13m-90}\leq 210$. Hence $(q,m)\in\Gamma_{a,b}(3,1)$ if $q^m>(3\cdot 2^{10})^{210}$, which clearly holds true whenever $\om(q^m-1)\geq 9632$. So we may consider the cases $\om(q^m-1)\leq 9631$. We now use Proposition \ref{P3.6} to reduce the lower bound mentioned above. For the sake of it, take $f=x^m-1$ so that $\tilde{f}$=gcd$\Bigg(f,\frac{x^m-1}{g}\Bigg)=\frac{x^m-1}{g}$ and $e$ is assumed to be the product of least $18$ prime divisors of $Q$, that is, $W(l)=2^{18}$. Then $s_1\leq 9613$ and $\delta$ assumes its minimum positive value when $\{p_{1},p_{2},\ldots,p_{9613}\}=\{67,71,\ldots,100483\}$, which gives $\delta>0.008195$ and $\Delta <1.18\times 10^6$. Thus, $3W(l)W(\tilde{f})$$\Delta<9.45\times 10^{14}=T$(say) and by Proposition \ref{P3.6}, it follows that $(q,m)\in\Gamma_{a,b}(3,1)$ if $q^{\frac{m}{2}-3}>T$, that is, if $q^m>T^{\frac{2m}{m-6}}$. Now, $m\geq 7$ implies that $\frac{2m}{m-6}\leq 14$. Hence, $(q,m)\in\Gamma_{a,b}(3,1)$ whenever $q^m>4.53\times 10^{209}$, which is true for $\om(q^m-1)\geq 97$. Repeat the process in a similar manner, it follows that $(q,m)\in\Gamma_{a,b}(3,1)$ if we have $q^{\frac{m}{2}-3}>2.75\times 10^7$. Thus $(q,m)\in\Gamma_{a,b}(3,1)$ unless $m=11$ and $q<945$; $m=10$ and $q<5244$; $m=9$ and $q<91107$; $m=8$ and $q<2.75\times 10^7$; $m=7$ and $q<7.56\times 10^{14}$. For each of these possible exceptions $(q,m)$ so that $m=9,10, 11$ and $3|q^m-1$, by testing Inequality $(\ref{E9})$ followed by $(\ref{E10})$ gives us that $(q,m)\in\Gamma_{a,b}(3,1)$ unless $(q,m)$ is equal to $(11,10)$. In addition, following same steps we have $(q,8)\in\Gamma_{a,b}(3,1)$ unless $q=$8, 11, 13, 16, 17, 19, 23, 25, 29, 32, 37, 41, 43, 47, 49, 53, 59, 83, 89.  Summarizing we get the following lemma.
	\begin{lem}
		Let $m=8,9,10,11$ and $q\geq 8$ be a prime power such that $3|q^m-1$. Then, $(q,m)\in\Gamma_{a,b}(3,1)$ unless the following possibilities.
		\begin{itemize}
			\item[1.] $m=8$ and $q=8, 11, 13, 16, 17, 19, 23, 25, 29, 32, 37, 41, 43, 47, 49, 53, 59, 83, 89$.
			\item[2.] $m=10$ and $q=11$.
		\end{itemize} 
	\end{lem}
	\begin{lem}
		Let $m=7$ and $q\geq 8$ be a prime power such that  $3|q^m-1$. Then, $(q,7)\in\Gamma_{a,b}(3,1)$ except possibly when $q\in \{13, 16, 19, 25, 31, 37, 43, 49, 61, 64, 67, 73,\\ 79, 97, 103, 109, 121, 127, 139, 151, 169, 181, 193, 199, 211, 223, 256, 277, 307, 337, 361,\\ 379, 421, 433, 463, 529, 547, 631, 673, 883\}$. 
	\end{lem}
	\bp
	Clearly, we have gcd$\Bigg(q-1,\frac{q^7-1}{q-1}\Bigg)=1$ or $7$. Further, any prime divisor of $\frac{q^7-1}{q-1}$ is either $7$ or of the form $7k+1$. We split the discussion into following two parts:\\
	\textbf{Case 1.} Let gcd$\Bigg(q-1,\frac{q^7-1}{q-1}\Bigg)=1$, which gives $Q=\frac{q^7-1}{q-1}$. Take $e=$gcd$(Q,7\cdot 23\cdot 43)$ and $f=1$ in Proposition \ref{P3.6} so that $s_2\leq 7$, $W(e)\leq 8$ and $W(\tilde{f})=1$. Prior to the above lemma, we have established that $(q,7)\in\Gamma_{a,b}(3,1)$ for $q>7.56\times 10^{14}$, which holds if $\om(q^7-1)\geq 56$. Note that $\om(Q)\geq 55$ implies that $\om(q^7-1)\geq 56$. So, we only deal with the cases $\om(Q)\leq 54$. Thus, any prime $p_i$ dividing $Q$ but not $e$ is of the form $7j+1$, for some $j\geq 7$.
    If we assume $q\geq 10^5$, then $\delta$ will assume it's minimum value, say $\delta_{min}$, whenever $\{p_1,p_2,\ldots,p_{s_1}\}$ be the set of consecutive $54$ primes of the form $7j+1$, $j\in\mathbb{N}$ starting with $p_1=71$. Thus $\delta>1-\sum_{i=1}^{54}\frac{1}{p_i}-(7/10^5)>0.902377=\delta_{min}$ and $\Delta<68.491$.
This gives us $(q,9)\in\Gamma_{a,b}(3,1)$, if $q>2.71\times 10^6$. Repeating the process in the similar way, we find that $(q,9)\in\Gamma_{a,b}(3,1)$ unless $q<780097$. Now, for each of the prime power $q$ satisfying $8\leq q< 780097$, first we test Inequality $(\ref{E9})$ and then for those not satisfy the same, we verify $(\ref{E10})$ and get the possible exceptions as $q=$13, 16, 19, 25, 31, 37, 49, 61, 67, 73, 79, 97, 103, 109, 121, 139, 151, 181, 193, 199, 223, 256, 277, 307, 361, 433, 529.\\
	\textbf{Case 2.} Let gcd$\Bigg(q-1,\frac{q^7-1}{q-1}\Bigg)=7$. Clearly, in this case, we have $7\nmid Q$ and divisors of $Q$ are of the form $7j+1$ for some $j\in\mathbb{N}$. Here, we choose $e=$gcd$(Q,23\cdot 43\cdot 71)$ and $f=1$. Following the same steps as mentioned in the previous case, we can get the possible exceptions as $q=43, 64, 127, 169, 211, 337, 379, 421, 463, 547, 631, 673, 883.$
    
    \ep
	\subsection{The cases $1\leq m\leq 4$:} From the definition of $k$-normal elements, we get that $0\leq k \leq m-1$. Clearly, it makes sense to discuss the existence of $3$-primitive $1$-normal polynomial for $m\geq 2$. We now verify the existence of a $3$-primitive $1$-normal element $\xi\in\Field_{q^{n}}$ such that $\mathrm{Tr}_{\Field_{q^{m}}/\Field_{q}}(\xi^{-1})=ab^{-1}$ and $\mathrm{N}_{\Field_{q^{m}}/\Field_{q}}(\xi)=b$, for any $a\in\Field_q$, $b\in\Field_q^*$, and for $m=2,3,4$. 
	\begin{lem}
		Let $m=2$ and $q\geq 2$ be a prime power such that $3|q^m-1$. Then $(q,2)\notin\Gamma_{a,b}(3,1)$.
	\end{lem}
	\bp
	Suppose that there exists a $3$-primitive $1$-normal element $\xi\in\Field_{q^{2}}$ such that $\mathrm{Tr}_{\Field_{q^{2}}/\Field_{q}}(\xi^{-1})=0$ and $\mathrm{N}_{\Field_{q^{2}}/\Field_{q}}(\xi)=b$, for some $b\in\Field_q^*$. This implies that $\xi^{2(q-1)}=1$, that is, ord$(\xi)\leq 2(q-1)<\frac{q^2-1}{3}$ for $q>5$ and thus $(q,2)\notin\Gamma_{a,b}(3,1)$ for $q> 5$. In $\Field_{4}$, $1$ is the only $3$-primitive element and $\mathrm{Tr}_{\Field_{4}/\Field_{2}}(1)=0$ and thus $(2,2)\notin\Gamma_{a^*,b}(3,1)$, for any $a^*\in\Field_q^*$. Note that when $m=2$, for any $1$-normal element $\xi\in\Field_{q^m}$ over $\Field_{q}$, we have either $\xi^q=\xi$ or $\xi^q=-\xi$. \\
	\textbf{Case 1.} Let $q=4$. Then, if $\xi\in\Field_{16}$ be any $3$-primitive $1$-normal element over $\Field_{4}$, then either $\xi^4=\xi$ or $\xi^4=-\xi$ and ord$(\xi)=5$. If $\xi^4=\xi$, then $\xi$ can not be $3$-primitive. Again, $\xi^4=-\xi$ and ord$(\xi)=5$ together implies that $\xi^4=1$, a contradiction. Hence, $(4,2)\notin\Gamma_{a,b}(3,1)$.\\ 
	\textbf{Case 2.} Let $q=5$. Then, if $\xi\in\Field_{25}$ be any $3$-primitive $1$-normal element over $\Field_{5}$, proceeding similarly as above, we get that $\xi^5=-\xi$ and ord$(\xi)=8$. Now, $\xi^5=-\xi$ gives us $\mathrm{Tr}_{\Field_{25}/\Field_{5}}(\xi^{-1})=0$, which holds for any $3$-primitive $1$-normal element over $\Field_{5}$. Hence, $(5,2)\notin\Gamma_{a,b}(3,1)$. 
	\ep
	\begin{lem}\lb{L4.9}
		Let $m=3$ and $q\geq 2$ be a prime power such that $3|q^m-1$. Then $(q,3)\notin\Gamma_{a,b}(3,1)$.
	\end{lem}
	\bp
	In this case, we must have $q\equiv 1(mod~ 3)$ and thus gcd $(q,3)=1$. This gives us, $x^3-1=(x-1)(x^2+x+1)$ or $x^3-1=(x-1)(x-c)(x-c^{-1})$ for some $c\notin\{0,1\}$. Let $\xi\in\Field_{q^3}$ an $3$-primitive $1$-normal element over $\Field_q$ such that $\mathrm{Tr}_{\Field_{q^{3}}/\Field_{q}}(\xi^{-1})=ab^{-1}$ and $\mathrm{N}_{\Field_{q^{3}}/\Field_{q}}(\xi)=b$, for some $a\in\Field_{q}$, $b\in\Field_q^*$. Since $\xi$ is $1$-normal, we must have deg(Ord$_q(\xi))= 2$.  \\
	\textbf{Case 1.} If Ord$_q(\xi)=x^2+x+1$, then we have $\mathrm{Tr}_{\Field_{q^{3}}/\Field_{q}}(\xi)=0$. Additionally,  $\mathrm{Tr}_{\Field_{q^{3}}/\Field_{q}}(\xi^{-1})=0$ gives us the minimal polynomial of $\xi\in\Field_{q^3}$ is of the form $x^3+d$ for some $d\in\Field_{q}$. This gives $\xi^{3(q-1)}=1$, that is, ord$(\xi)\leq 3(q-1)< \frac{q^3-1}{3}$ for $q>2$, a contradiction. \\
	\textbf{Case 2.} If Ord$_q(\xi)=(x-1)(x-c)$, then we have $(\xi^{q}-c\xi)^{q-1}=1$, that is, $\xi^{q}-c\xi\in\Field_{q}$. Let $\xi^{q}=c\xi+f$, for some $f\in\Field_{q}$. Now, $\mathrm{Tr}_{\Field_{q^{3}}/\Field_{q}}(\xi^{-1})=0$ gives us the minimal polynomial of $\xi\in\Field_{q^3}$ is of the form $x^3+dx^2+e$ for some $d,e\in\Field_{q}$. Then, $\xi^{3}=d_{1}{\xi}^2+e_1$ for some $d_1,e_1\in\Field_{q}$, which implies that $\xi^{3q}=d_{1}{\xi}^{2q}+e_1$, that is, $(c\xi+f)^3=d_1(c\xi+f)^2+e_1$. Further, by substituting $\xi^{3}=d_{1}{\xi}^2+e_1$, we get a polynomial $P$ of degree $\leq 2$ over $\Field_{q}$ such that $P(\xi)=0$. Then $\xi\in\Field_{q^{2}}$, that is, ord$(\xi)\leq q^2-1<\frac{q^3-1}{3}$ for $q>3$, a contradiction.\ep
	\begin{lem}
		If there exists a $3$-primitive $1$-normal element $\xi\in\Field_{q^{4}}$ over $\Field_{q}$ such that $\mathrm{Tr}_{\Field_{q^{4}}/\Field_{q}}(\xi^{-1})=ab^{-1}$ and $\mathrm{N}_{\Field_{q^{4}}/\Field_{q}}(\xi)=b$, for some $a\in\Field_{q}$, $b\in\Field_q^*$, then  $3|q^4-1$ and $q\equiv 1 (mod 4)$.
	\end{lem}
	\bp
	Let $\xi\in\Field_{q^4}$ an $3$-primitive $1$-normal element over $\Field_q$ such that $\mathrm{Tr}_{\Field_{q^{4}}/\Field_{q}}(\xi^{-1})=ab^{-1}$ and $\mathrm{N}_{\Field_{q^{4}}/\Field_{q}}(\xi)=b$, for some $a\in\Field_{q}$, $b\in\Field_q^*$. If $q\not\equiv 1(mod 4)$, then $x^4-1=(x-1)^4$ or $(x-1)(x+1)(x^2+1)$. Since $\xi$ is $1$-normal, we must have deg(Ord$_q(\xi))= 3$. We now consider the following possible cases:\\ 
	\textbf{Case $1$.} If Ord$_q(\xi)=x^3+x^2+x+1$, then we have $\mathrm{Tr}_{\Field_{q^{4}}/\Field_{q}}(\xi)=0$. In addition to this, $\mathrm{Tr}_{\Field_{q^{4}}/\Field_{q}}(\xi^{-1})=0$ gives us the minimal polynomial of $\xi$ is of the form $x^4=cx^2+d$, for some $c,d\in\Field_{q}$.   Then ${\xi}^4=c{\xi}^2+d$, which implies that ${\xi}^{4q}=c{\xi}^{2q}+d$, which further implies that $\xi^{4q}-\xi^4=c\cdot(\xi^{2q}-\xi^2)$. Since $\xi^{2q}\neq \xi^2$, it follows that $\xi^{2q}+\xi^2=c$, which gives $\xi^{2q^2}+\xi^{2q}=c$. Thus $\xi^{2(q^2-1)}=1$, that is, ord$(\xi)\leq 2(q^2-1)<\frac{q^4-1}{3}$ for $q>2$.\\
	\textbf{Case $2$.} If Ord$_q(\xi)=(x-1)(x^2+1)$, then we get that $(\xi^{q^2}+\xi)^{q-1}=1$, that is, $\xi^{q^2}+\xi=g$ for some $g\in\Field_{q}$. Since $\xi$ is $3$-primitive, we must have $g\neq 0$. Further, $\mathrm{Tr}_{\Field_{q^{4}}/\Field_{q}}(\xi^{-1})=0$ gives us the minimal polynomial of $\xi$ is of the form $x^4-dx^3-ex^2-f$ for some $d,e,f\in\Field_{q}$. Then, we have ${\xi}^4=d{\xi}^3+e{\xi}^2+f$, which gives that ${\xi}^{4q^2}=d{\xi}^{3q^2}+e{\xi}^{2q^2}+f$. Thus $({\xi}^{4q^2}-{\xi}^4)=d({\xi}^{3q^2}-{\xi}^3)+e({\xi}^{2q^2}-{\xi}^2)$, which implies that $(\xi^{q^2}+\xi)(\xi^{2q^2}+\xi^2)=d(\xi^{2q^2}+\xi^{q^2+1}+\xi^2)+e(\xi^{q^2}+\xi)$. Further, by substituting $\xi^{q^2}=g-\xi$, we get a polynomial $R$ of degree $2$ over $\Field_{q}$ such that $R(\xi)=0$. Hence ord$(\xi)\leq q^2-1<\frac{q^4-1}{3}$ for $q\geq 2$, a contradiction.\\

	Thus, a $3$-primitive $1$-normal element $\xi\in\Field_{q^4}$ with the properties as mentioned above may exist, only if $q\equiv 1 (mod ~4)$ and $3|q^4-1$.    
	\ep
\section{Open Problems}
We hereby conclude the article by stating the following open problems:
\begin{itemize}
\item  For $m=4$, are the conditions $q^m\equiv 1(mod ~3)$ and $q\equiv 1 (mod ~4)$ sufficient to ensure the existence of a $3$-primitive $1$-normal polynomial over $\Field_q$ of degree $m$ with last two prescribed coefficients ${?}$ \\
\item For $m=5,6$, does there exists a $3$-primitive $1$-normal polynomial over $\Field_q$ of degree $m$ with last two prescribed coefficients ${?}$ 
\end{itemize}
	
	\section{Statements and Declarations}
	There are no known financial interests or personal relationships that could have influenced the findings of this paper.
	\section{Acknowledgments}
 First author is supported by the National Board for Higher Mathematics (NBHM), Department of Atomic Energy (DAE), Government of India, under Ref No. 0203/6/\\2020-R\&D-II/7387. All authors are equally contributed.
	
\end{document}